\documentclass{article}

\usepackage{amsmath} 
\usepackage{amssymb}  
\usepackage{amsthm}
\usepackage{latexsym}
\usepackage{amsfonts}
\usepackage{xcolor}
\usepackage{graphicx}
\usepackage{mathtools}
\usepackage{amsfonts,amsmath,amssymb,mathrsfs}

\usepackage[letterpaper,top=2cm,bottom=2cm,left=3cm,right=3cm,marginparwidth=1.75cm]{geometry}
\usepackage[colorlinks=true, allcolors=blue]{hyperref}

\newtheorem{thm}{Theorem}

\newtheorem{prop}[thm]{Proposition}
\theoremstyle{definition}
\newtheorem{defn}[thm]{Definition}

\newtheorem{ex}[thm]{Example}

\newcommand{\calB}{\mathcal B} 
\newcommand{\calJ}{\mathcal J}
\newcommand{\calH}{\mathcal H} 
\newcommand{\calD}{\mathcal D}
\newcommand{\calG}{\mathcal G}
\newcommand{\calR}{\mathcal R}

\newcommand{\R}{\mathbb R}
\newcommand{\N}{\mathbb N}
\newcommand{\tr}{\operatorname{tr}}

\newcommand{\pr}{\operatorname{pr}^{(m_d)}}
\newcommand{\pz}{\partial_z}
\newcommand{\calI}{\mathcal I}
\newcommand{\J}{\mathbb J}

\title{Jet space extensions of infinite-dimensional Hamiltonian systems\footnote{This work was initiated during a research stay of TP at the LAGEPP laboratory at UCBL Lyon 1. TP acknowledges support by the European Union by means of an Erasmus+ scholarship.}} 
\author{T. Preuster$^1$ and M. Schaller$^1$ and B. Maschke$^2$}
\date{$^1$Institute of Mathematics, Technische Universit\"at Ilmenau, Germany\\
$^2$Universit{\'e} Claude Bernard Lyon 1, CNRS, LAGEPP UMR 
	5007, France\\[2ex]%
    \today}
\begin{document}
\maketitle

\begin{abstract}
\noindent     We analyze infinite-dimensional Hamiltonian systems corresponding to partial differential equations on one-dimensional spatial domains formulated with formally skew-adjoint Hamiltonian operators and nonlinear Hamiltonian density. 
    In various applications, the Hamiltonian density can depend on spatial derivatives of the state 
    such that these systems can not straightforwardly be formulated as boundary port-Hamiltonian system using a Stokes-Dirac structure. 
    In this work, we show that any Hamiltonian system of the above class can be reformulated as a Hamiltonian system on the jet space, in which the Hamiltonian density only depends on the extended state variable itself and not on its derivatives. Consequently, well-known geometric formulations with Stokes-Dirac structures are applicable. 
    Additionally, we provide a similar result for dissipative systems. We illustrate the developed theory by means of the the Boussinesq equation, the dynamics of an elastic rod and the Allen-Cahn equation.
\end{abstract}


\section{Introduction}\label{sec:intro}
Hamiltonian systems encapsulate the evolution of physical systems in terms of their energies and associated generalized coordinates. At their core lies the Hamiltonian function, often described as the total energy storage of a system. Based on Hamiltonian systems, the theory of port-Hamiltonian systems was developed to open this class to external influence, \cite{maschke1992intrinsic, maschke1993port}. 
The theory has developed into a broad and active field of research and port-Hamiltonian systems have also been considered on infinite dimensional spaces \cite{le2005dirac, JacoZwar12}. 

As for finite dimensional systems, the central aspects of infinite-dimensional port-Hamiltonian systems are energy storage and energy routing. Energy routing can be modeled using Stokes-Dirac structures \cite{le2005dirac}, whereas energy storage is captured via the Hamiltonian function. Frequently, when modeling distributed parameter systems, the storage function involves derivatives of the state variables. This means that the Hamiltonian no longer solely depends on the state variables; rather, it also incorporates spatial derivatives of the state.

It has been observed in several applications that the Hamiltonian form is preserved when the system is embedded in the jet space, \cite{maschke2013alternative}. In addition, it was even possible to represent it as a boundary port-Hamiltonian system, \cite{MascvdSc23}. In this work, we prove this property for conservative and dissipative port-Hamiltonian systems. 
Port-Hamiltonian systems have also been formulated directly on jet bundles~\cite{schoberl2014jet}, where the authors allow for multi-dimensional spatial domains.

This article is organized as follows. In Section \ref{sec:problem}, we formulate the problem statement followed by motivating examples from fluid dynamics and mechanics. In Section \ref{sec:jet} we provide the first main result of this work, showing that jet space extensions of Hamiltonian systems are again Hamiltonian, and the extended Hamiltonian density does not depend on derivatives of the state. Further, in Section \ref{sec:diss}, we extend this result to dissipative Hamiltonian systems. Last, in Section \ref{sec:outlook}, we conclude this work and provide an outlook to future research perspectives.\medskip

\textbf{Notation.}
We define the $m_d$-th prolongation of a state variable $x$ by $\pr x \coloneqq (x, \partial_z x, \ldots, \partial_z^{m_d} x)$. We will write $[n:m] \coloneqq [n,m]\cap \mathbb{Z}$ for appropriate index sets. For a matrix differential operator $\calJ$ as in \eqref{eq:Hamiltonian_operator} we denote the 
$(i,j)-$th entry by  $\calJ_{ij}$ and by $\calJ_{i*}$ the $i$-th row. 
We consider index sets  $\calI \coloneqq \{(i_1, j_1), \ldots, (i_l, j_l)\} \subset [1:n] \times [0:m_d]$ for $q,m_d \in \N$ and denote by $\operatorname{ind}: \{j_1, \ldots, j_l\} \to [1:l]$ the choice function mapping any $j_k$ to $\operatorname{ind}(j_k) \coloneqq k$ for $k \in [1:l]$. Further, we let $i: \{j_1, \ldots, j_l\} \to \{i_1, \ldots, i_l\}$ defined by $j_k \mapsto i_k$ for $k \in [1:l]$. 

\section{Motivation and problem statement}\label{sec:problem} 

In this section, we first recall the definition of Hamiltonian systems on a one-dimensional spatial domain where the Hamiltonian operator is a constant coefficient matrix operator and the Hamiltonian functional
depends on the spatial derivatives of the state variables, i.e. it is defined on the prolongation of the state space \cite{olver}. We then briefly recall the definition of boundary port-Hamiltonian
systems \cite{van2002hamiltonian,le2005dirac} with respect to a Stokes-Dirac structure associated with the Hamiltonian operator, which are defined when the Hamiltonian functional depends only on the spatial
variable and the state variables and not its derivatives. Finally, we recall on two simple examples, how one may use the jet space extension of a Hamiltonian system in order to formulate a boundary port-Hamiltonian
systems, as a motivation to derive a general procedure.

\subsection{Hamiltonian systems}\label{subsec:Hamiltonian-systems}

Let an interval $[a,b]\subset\R$ be given. In this paper, we consider the matrix differential operators of the form
\begin{align}\label{eq:Hamiltonian_operator}
\calJ=\sum_{j=0}^{m}P_{j}\partial_{z}^{j},
\end{align}
for $P_{j}\in\R^{n\times n}$ with $P_{j}=(-1)^{j+1}P_{j}^{\top}$ for $j\in[0:m]$. Using integration by parts,  $\calJ$ is formally skew-adjoint in the $L^2([a,b];\R)$, cf.~\cite{MascvdSc23}.
\begin{defn}[Hamiltonian systems]
    Let $\calJ$ be as in \eqref{eq:Hamiltonian_operator} and let $H:(\R^n)^{m_d} \to \R$
    be a Hamiltonian density depending on the state $x:[a,b]\to\R$ and its first $m_{d} \in \N$ spatial derivatives, denoted by the prolongation $\pr x$. An \textit{infinite-dimensional Hamiltonian system} generated by the Hamiltonian functional 
    \begin{equation}\label{eq:HamFunctional}
        \calH(\pr x)\coloneqq\int_{a}^{b}H(\pr x(z))\,\mathrm{d}z
    \end{equation}
     is given by the system dynamics 
    \begin{equation}\label{eq:ham_sys}
        \frac{\mathrm{d}}{\mathrm{d}t}x=\calJ\delta_{x}\calH(\pr x),
    \end{equation}
    where $\delta_{x}\calH(\pr x)$ denotes the \textit{variational derivative} of $\calH(\pr x)$ w.r.t.~$x$ defined by 
    \begin{align}\label{eq:VariationalDerivative}
      \left.\frac{d}{d\varepsilon}\right|_{\varepsilon=0}\calH[x+\varepsilon\eta]=\int_a^b \langle \delta_{x}\calH(\pr x)(z), \eta(z) \rangle_{\R^n} \, \mathrm{d}z
    \end{align}
    for all $\eta \in  H_{0}^{m_{d}}([a,b],\R^{n})$ \cite[Def. 4.1]{olver}.
\end{defn}
We proceed with two examples. 
\begin{ex}\label{ex:Elastic-rod}
We consider the dynamics of a one-dimensional elastic rod given by
\begin{align*}
    \frac{\mathrm{d}}{\mathrm{d}t}
    \begin{pmatrix}
        u\\ p
    \end{pmatrix}=
    \begin{pmatrix}
        0 & 1 \\ -1 & 0
    \end{pmatrix}
    \begin{pmatrix}
        \delta_{u} \calH(u,p,\partial_{z}u) \\ \delta_{p} \calH(u,p,\partial_{z}u)
    \end{pmatrix}
\end{align*}
with the displacement variable $u$ and the momentum density $p$ of the elastic rod \cite{MascvdSc23}. The Hamiltonian functional is defined by 
\begin{align*}
    \calH(u,p,\partial_{z}u) = \frac{1}{2}\int_{a}^{b}\frac{1}{(\rho A)}p(z)^{2}+k\,u(z)+T\left(\partial_{z}u(z)\right)^{2}\, \mathrm{d}z
\end{align*}
where the first term corresponds to the kinetic energy, the second is the the gravity potential and third is structural elastic energy. The parameters are elasticity modulus $T>0$, volume mass density $\rho$, section area $A>0$, and $k\geq0$.
\end{ex}

\begin{ex}\label{ex:Boussinesq-equation}
The Boussinesq equation 
\begin{align*}
    \frac{\mathrm{d}^{2}}{\mathrm{d}t^{2}}u=\frac{1}{3}\partial_{z}^{4}u+\frac{4}{3}\partial_{z}^{2}(u^{2})
\end{align*}
is a nonlinear approximation equation for water waves in shallow water whereas $u$ denotes the displacement of the surface.  The model can be rewritten as the following Hamiltonian system of two conservation laws \cite[p.~459 ff]{olver}
\begin{align}\label{eq:Boussinesq}
    \frac{\mathrm{d}}{\mathrm{d}t}
    \begin{pmatrix}
        u\\ v
    \end{pmatrix}=
    \begin{pmatrix}
        0 & \partial_{z}\\
        \partial_{z} & 0
    \end{pmatrix}
    \begin{pmatrix}
        \delta_{u}\calH(u,v,\partial_{z}u) \\ \delta_{v}\calH(u,v,\partial_{z}u)
    \end{pmatrix}
\end{align}
with the Hamiltonian functional 
\begin{align*}
    \calH(u,v,\partial_{z}u)=\int_{a}^{b}-\frac{1}{6}(\partial_{z}u(z))^{2}+\frac{4}{9}u(z)^{3}+\frac{1}{2}v(z)^{2}\,\mathrm{d}z.
\end{align*}
Here, whereas the Hamiltonian operator in \eqref{eq:Boussinesq} is a linear matrix operator, the energy functional is not quadratic as it involves a third order term in $u$. Further, note that it is indefinite, due to the negative sign of the first term. 
\end{ex}
\subsection{Boundary port-Hamiltonian systems}
Boundary Port-Hamiltonian systems are an extension of Hamiltonian systems \eqref{eq:ham_sys} when considering functions which do not vanish at the boundary and the operator $\calJ$ is no more skew-adjoint, thereby enabling the description of \emph{open} Hamiltonian systems which are interacting with their environment \cite{van2002hamiltonian}. Therefore the Hamiltonian operator $\calJ$ in \eqref{eq:Hamiltonian_operator} is extended to a Dirac structure\footnote{Let $\mathcal{F}$ and $\mathcal{E}$ be two real vector spaces endowed with a non-degenerate bilinear form $\left\langle \cdot,\cdot\right\rangle $, called \emph{pairing}.
On the product space $\mathcal{B}=\mathcal{F}\times\mathcal{E}$, called bond space, consider the \emph{plus pairing} $\ll(f_{1},e_{1}),(f_{2},e_{2})\gg_{+} \coloneqq \langle e_{1},\,f_{2}\rangle +\langle e_{2},\,f_{1}\rangle $. A \emph{Dirac structure} is a subspace $\mathcal{D}\subset\mathcal{B}$
such that $\mathcal{D}= \mathcal{D}^{\bot}$, with $\perp$ denoting
the orthogonal complement with respect to the plus pairing $\ll,\gg_{+}$.} \cite{Courant88, courant1990dirac} called \emph{Stokes-Dirac structure}.
\begin{prop}\label{prop:stokes_dirac}\cite{le2005dirac,villegas2007port,MascvdSc23} 
Consider the flow space $\mathscr{F}=L^{2}\left( \left[a,\,b\right],\,\R^n\right)\times \mathbb{R}^{p} \ni \left(f,\,f_{\partial}\right)^{\top}=\bar{f}$ and the effort space $\mathscr{E}=\mathscr{F}^*\ni\left(e ,\, e_{\partial}\right)^{\top}=\bar{e}$. The bond space is defined by $\mathscr{B}=\mathscr{F}\times\mathscr{E}$ and equipped with the symmetric pairing 
\begin{align*}
    \left\langle \bar{f},\,\bar{e}\right\rangle_{\mathscr{B}}
    =
    \left\langle 
    \begin{pmatrix}
        f\\f_{\partial}
    \end{pmatrix},
    \begin{pmatrix}
        e\\ e_{\partial}
    \end{pmatrix}\right\rangle_{\mathscr{B}}
    =\left\langle f,e\right\rangle _{L^{2}} - \left\langle f_{\partial},e_{\partial}\right\rangle_{\R^p}
\end{align*}
Consider a Hamiltonian operator $\mathcal{J}$ defined in \eqref{eq:Hamiltonian_operator}. Then there exists an integer $p$, $W \in \R^{2p\times2p}$ and a $(p\times n)$ matrix differential operator $\mathcal{T}$ such that the subspace 
\begin{equation}\label{eq:StokesDiracConstEq}
    \mathscr{D}_{\mathcal{J}} \hspace{-0.5ex} = \hspace{-0.5ex} \left\{ \hspace{-0.75ex}\left( \hspace{-0.95ex}
    \begin{pmatrix}
        f \\ f_{\partial}
    \end{pmatrix} \hspace{-0.5ex} ,\hspace{-0.5ex}   
    \begin{pmatrix}
        e \\ e_{\partial}
    \end{pmatrix}  \hspace{-0.95ex} \right) \hspace{-0.5ex} \in \mathscr{B} \middle| \hspace{-0.5ex}
    \begin{array}{l}
        e \in H^m([a,b],\R^n) \ \text{s.t.} \\ 
     f=\mathcal{J}e,
    \begin{pmatrix}
        f_\partial \\ e_{\partial}
    \end{pmatrix}\hspace{-0.5ex} = \hspace{-0.5ex} W \mathrm{tr} \left(\mathcal{T} e\right)
    \end{array} \hspace{-1ex} \right\}
\end{equation}
is a Dirac structure, called \emph{Stokes-Dirac structure}, with respect to a symmetric pairing on the bond space, i.e., $\ll\left(\bar{f}_{1},\,\bar{e}_{2}\right),\,\left(\bar{f}_{2},\,\bar{e}_{1}\right)\gg_{+}=\left\langle \bar{f}_{1},\bar{e}_{2}\right\rangle _{\mathscr{B}}+\left\langle \bar{f}_{2},\bar{e}_{1}\right\rangle _{\mathscr{B}}$ for $\left(\bar{f}_{i},\,\bar{e}_{i}\right)\in\mathscr{B}, i\in\left\{ 1,\,2\right\} $.
\end{prop}
\noindent The pairs $\left(f_{\partial},e_{\partial}\right)$ are called port variables and have various equivalent definitions, cf.~\cite{le2005dirac, villegas2007port, JacoZwar12, MascvdSc23}. Note that, by the definition of a Dirac structure, a consequence of
the Proposition \ref{prop:stokes_dirac} is that 
\begin{equation}\label{eq:generalizedSkewAdjointness}
    \int_{a}^{b}e^{T}\calJ e\,\mathrm{d}z+\int_{a}^{b}\left(\calJ e\right)^{T}e\,\mathrm{d}z-e_{\partial}^{\top}f_{\partial}=0
\end{equation}
for $(f,f_\partial, e, e_\partial) \in \mathscr{D}_\calJ$ which is indeed a generalization of the skew-adjointness of the operator
$\mathcal{J}$. Using the Stokes-Dirac structure, Hamiltonian systems \eqref{eq:ham_sys} have been extended to boundary port-Hamiltonian systems with respect to the Stokes-Dirac structure \eqref{eq:StokesDiracConstEq}, in the case when the Hamiltonian density $H(x)$ depends only on the state variables and not their spatial derivatives.

\begin{defn}\label{def:Linear_BPHS} [Boundary Port-Hamiltonian System \cite{schaftGeomPhys02} ] The  boundary port-Hamilto\-nian system
on the state space $\mathscr{X}=L_{2}\left([a,\,b],\,\mathbb{R}^{n}\right)$ with respect to the Stokes-Dirac structure $\mathscr{D}_{\mathcal{J}}$ defined in \eqref{eq:StokesDiracConstEq} associated with the Hamiltonian matrix operator $\mathcal{J}$ and Hamiltonian $\calH$ with density $H =H(x)$ is the system 
\begin{equation}\label{eq:BPHS_GeometricDef}
    \left(\left(
    \begin{array}{c}
        \tfrac{\mathrm{d} x}{\mathrm{d} t}\\
        f_{\partial}
    \end{array}\right),\left(
    \begin{array}{c}
        \delta_{x}\calH\\
        e_{\partial}
    \end{array}\right)\right) \in \mathscr{D}_{\mathcal{J}},
\end{equation}
where $\left(f_{\partial},\,e_{\partial}\right)$ are the \emph{boundary
port variables}.
\end{defn}
The definition of the Stokes-Dirac structure and the property \eqref{eq:generalizedSkewAdjointness}
yields the energy balance equation 
\begin{align}\label{eq:EnergyBalance_BPHS}
    \frac{\mathrm{d}}{\mathrm{d}t}\calH=e_{\partial}^{\top}f_{\partial}
\end{align}
If the Hamiltonian functional is bounded from below, the boundary port-Hamiltonian system is lossless.
\subsection{Problem statement}
However the very classic examples of the elastic rod and the Boussinesq equation, presented in the section \ref{subsec:Hamiltonian-systems}, have a Hamiltonian functional which does depend on derivatives of the state such that Definition~\ref{def:Linear_BPHS} may not be readily applied. Of course, if the Hamiltonian depends on the prolonged state, i.e., $\calH(\pr x)$, one may still define the dynamics according to \eqref{eq:BPHS_GeometricDef} but the energy balance equation~\eqref{eq:EnergyBalance_BPHS} is no more satisfied \cite{MascvdSc23} and some interface variables are missing to render the system complete. 
\begin{ex}
    Consider the Hamiltonian formulation of the Korteweg–De Vries equation
    \begin{equation*}
        \partial_t x = \partial_z \delta_x \calH_1(x, \partial_z x),  
    \end{equation*}
    with $\calH_1(x, \partial_z x):=\int_a^b -\frac{1}{2}(\partial_z x)^2 + \frac{1}{6}x^3 \, \mathrm{d}z$, \cite{maschke2013alternative}. The Hamiltonian operator $\partial_z$ induces a Stokes-Dirac structure but one may compute
    \begin{align*}
           \frac{\mathrm{d}}{\mathrm{d}t}\calH=e_{\partial}^{\top}f_{\partial} + \left[\partial_z x \delta_x \calH(x, \partial_z x) \right]_a^b
    \end{align*}
    for the corresponding boundary port variables $f_\partial, e_\partial \in \R$.
\end{ex}
There are two ways to deal with this problem. The first one is to introduce additional port variables associated with the Hamiltonian functional $\calH(\pr x)$ by defining Stokes-Lagrange subspaces and a boundary port-Hamiltonian systems with two types of port variables: power ports associated with the Stokes-Dirac structure and energy ports associated with the Hamiltonian $\calH(\pr x)$~\cite{MascvdSc23}. 
The second way which we use in this paper, is to lift the system to a higher dimensional state including the spatial derivatives of the state appearing in the Hamiltonian functional. This lifting is commonly called a \textit{jet space extension}, \cite[p. 96]{olver}. This has been actually utilized in the original approach for defining boundary port-Hamiltonian systems for mechanical systems \cite{schaftGeomPhys02,le2005dirac} or fluid-mechanical systems \cite{maschke2013alternative}.

We will now illustrate this lifting procedure by means of the two applications of Example~\ref{ex:Elastic-rod} and \ref{ex:Boussinesq-equation}.
\begin{ex}\label{ex:ElasticRod_Prolongation} Let us consider again the example of the flexible rod presented in Example \ref{ex:Elastic-rod} where the Hamiltonian, the total elasto-dynamic energy, depends on the first derivative of the displacement $\partial_{z}u$. Define the strain variable as $\varepsilon(t,z)\coloneqq\partial_{z}u(t,z)$ and consider the augmented state 
\begin{align*}
    \bar{x}=
    \begin{pmatrix}
        \bar{x}_{1} &
        \bar{x}_{2} &
        \bar{x}_{3}
    \end{pmatrix}^T \coloneqq
    \begin{pmatrix}
        u&
        p&
        \varepsilon
    \end{pmatrix}^T
\end{align*}
with the Hamiltonian functional defined in the lifted variables 
$$\bar{\calH}(\bar{x})=\frac{1}{2}\int_{a}^{b}\frac{1}{(\rho A)}\bar{x}_{2}(z)^{2}+k\bar{x}_{1}(z)^{2}+T\left(\bar{x}_{3}(z)\right)^{2}\,\mathrm{d}z.$$
Following \cite{MascvdSc23} we can express the dynamics of the elastic rod equivalently by 
\begin{align*}
    \frac{\mathrm{d}}{\mathrm{d}t}\bar{x}=
    \begin{pmatrix}
        0 & 1 & 0\\
        -1 & 0 & \partial_{z}\\
        0 & \partial_{z} & 0
    \end{pmatrix}\delta_{\bar{x}}\bar{\calH}(\bar{x}).
\end{align*}
Again, we obtain a Hamiltonian system on the jet space in which, by construction, the Hamiltonian only depends on the state and not on its derivatives. \end{ex}

\begin{ex}\label{ex:Boussinesq-equation_Prolongation} Consider now Example \ref{ex:Boussinesq-equation} of the Boussinesq equation where the Hamiltonian functional depends also on $\partial_{z}u$. Thus, define the augmented state $\bar{x}$ as in \cite{maschke2013alternative}
by 
\begin{align}\label{eq:Boussinesq_state}
    \bar{x}=
    \begin{pmatrix}
        \bar{x}_{1}&
        \bar{x}_{2}&
        \bar{x}_{3}
    \end{pmatrix}^T\coloneqq
    \begin{pmatrix}
        u&
        v&
        \pz u
    \end{pmatrix}^T
\end{align}
and the Hamiltonian functional by
\begin{align}\label{eq:Boussinesq_Hamiltonian}
    \bar{\calH}(\bar{x})=\int_{a}^{b}-\frac{1}{6}(\bar{x}_{3}(z))^{2}+\frac{4}{9}\bar{x}_{1}(z)^{3}+\frac{1}{2}\bar{x}_{2}(z)^{2}\,\mathrm{d}z.
\end{align}
The dynamics \eqref{eq:Boussinesq} on the augmented state space may then again be written as a Hamiltonian system 
\begin{align}\label{eq:jet_Boussineq}
    \frac{\mathrm{d}}{\mathrm{d}t}\bar{x}=
    \begin{pmatrix}
        0 & \partial_{z} & 0\\
        \partial_{z} & 0 & -\partial_{z}^{2}\\
        0 & \partial_{z}^{2} & 0
    \end{pmatrix}\delta_{\bar{x}}\bar{\calH}(\bar{x})
\end{align}
in which the Hamiltonian $\bar{\calH}$ does only depend on the extended state $\bar{x}$ and not on its spatial derivatives. \end{ex}

In the sequel, we shall give a constructive method for the prolongation of the Hamiltonian systems \eqref{eq:ham_sys} to Hamiltonian systems on the jet space. 

\section{Hamiltonian jet space extensions}\label{sec:jet}
In this section, we give a constructive method to derive the lift of Hamiltonian system to the jet space defined by the arguments of the Hamiltonian function.

\subsection{Definition of the extended state space}
In the first instance, we shall define the extended state space on which the Hamiltonian system will be lifted. Therefore, we consider a subset of indices 
\begin{align*}
    \calI\coloneqq\{(i_{1},j_{1}),\ldots,(i_{l},j_{l})\}\subset[1:n]\times[0:m_{d}]
\end{align*}
and assume that $H$ depends only on the variable set $V_{\calI}\coloneqq\{\pz^{j}x_{i}\,|\,(i,j)\in\calI\}$. 
Further, we vectorize $V_{\calI}$ consistently to the notation in~\eqref{eq:ham_sys} by first ordering by the order of differentiation and second by coordinate. By definition of this index set, we write $H(\pr x)=H(V_{\calI})$. We first require a result on the variational derivative of this Hamiltonian functional.
\begin{prop}\label{prop:vario_deriv} 
Let $\calH(V_{\calI})=\int_{a}^{b}H(V_{\calI})\,\mathrm{d}z$. Then 
\begin{align*}
    \delta_{x}\calH(V_{\calI})=\sum_{j=1}^{m_{d}}(-\pz)^{j}\frac{dH}{d(\pz^{j}x)}(V_{\calI}).
\end{align*}
\end{prop}
\begin{proof} The claim follows directly from the definition of the variational derivative \eqref{eq:VariationalDerivative} together with the chain rule and integration by parts:
\begin{equation*}\label{eq:var_deriv_calc}
    \delta_x \calH (V_{\calI}) = \left(\begin{matrix}
        \sum\limits_{j=1}^{m_d} (-\pz)^j \frac{d H}{d (\pz ^j x_{1})} \\ \vdots \\ \sum\limits_{j={1}}^{m_d} (-\pz)^{j} \frac{d H}{d (\pz ^j x_{n})} 
    \end{matrix}\right)
     = \sum_{j=1}^{m_d} (-\pz)^j \left(\begin{matrix}
          \frac{d H}{d (\pz ^j x_{1})} \\ \vdots \\  \frac{d H}{d (\pz ^j x_{n})}
     \end{matrix}\right).
\end{equation*} 
\end{proof} 
Now define the augmented state 
\begin{equation*}
    \bar{x}
   \coloneqq \begin{pmatrix}
        \partial_z^{j_1} x_{i_1} & \partial_z^{j_2} x_{i_2} & \ldots & \partial_z^{j_l} x_{i_l}
    \end{pmatrix}^\top.
\end{equation*}
Having this substitution in mind, we can rewrite the Hamiltonian as a function in the lifted state 
\begin{align}\label{eq:lifted_Ham}
    \bar{\calH} (\bar{x}) \coloneqq \int_a^b H(\bar{x}_1
         , \bar{x}_2,  \ldots  \bar{x}_{l}) \, \mathrm{d}z.
\end{align}
It is important to note that from now on we interpret the Hamiltonian as a function in several variables without dependencies between the functions $\bar{x}_1, \ldots  \bar{x}_{l}$.

\subsection{Lift of the Hamiltonian system}

We are now in the position to formulate the central result of this note, showing how to lift a Hamiltonian system to a higher dimensional Hamiltonian system on the jet space $V_{\calI}$. We set  $D_\pm \coloneqq \operatorname{diag}( (\pm \pz)^{j_1},(\pm \pz)^{j_2}, \ldots, (\pm \pz)^{j_l})$.

\begin{prop}\label{prop:Extension1} 
The extension of \eqref{eq:ham_sys} onto the jet space is again a Hamiltonian system, given by 
\begin{equation}\label{eq: explicit_jet_dynamics}
    \frac{\mathrm{d}}{\mathrm{d}t}  
    \bar{x}
    =\mathbb{J} \delta_{\bar{x}} \bar{\calH}(\bar{x}) 
\end{equation}
with the formally skew-adjoint differential operator 
\begin{equation}\label{eq:ExtendedHamiltonianOperator}
     \mathbb{J} = 
     \begin{pmatrix}
        \pz^{j_1} \calJ_{i_1 i_1}(-\pz)^{j_1} & \ldots & \pz^{j_1} \calJ_{i_1 i_l}(-\pz)^{j_l} \\  
        \vdots & \ddots & \vdots \\ 
        \pz^{j_l} \calJ_{i_l i_1}(-\pz)^{j_1} & \ldots  & \pz^{j_l} \calJ_{i_l i_l}(-\pz)^{j_l}
    \end{pmatrix}.
\end{equation}
\end{prop}
\begin{proof} 
We consider the system dynamics \eqref{eq:ham_sys}
which extends to the dynamics of the states $V_I$:
\begin{equation*}
    \frac{\mathrm{d}}{\mathrm{d}t} \hspace{-1mm}\left( \begin{matrix}
        \partial_z^{j_1} x_{i_1} \\ \partial_z^{j_2} x_{i_2} \\ \vdots \\ \partial_z^{j_l} x_{i_l}
    \end{matrix}\right) \hspace{-1mm}
    = D_+\frac{\mathrm{d}}{\mathrm{d}t} \hspace{-1mm} \left( \begin{matrix}
          x_{i_1} \\  x_{i_2} \\ \vdots \\ x_{i_l}
    \end{matrix}\right)\hspace{-1mm}=D_+ \hspace{-1mm}\left(\begin{matrix}
                \calJ_{i_1*} \\
                \calJ_{i_2*} \\
                \vdots \\
                \calJ_{i_l*} 
            \end{matrix}\right) \delta_x \calH (V_{\calI}),
\end{equation*}
by Schwarz' theorem. Then Proposition \ref{prop:vario_deriv} yields that
\begin{align*}
     \frac{\mathrm{d}}{\mathrm{d}t} \hspace{-1mm}\left( \begin{matrix}
        \partial_z^{j_1} x_{i_1} \\ \partial_z^{j_2} x_{i_2} \\ \vdots \\ \partial_z^{j_l} x_{i_l}
    \end{matrix}\right) \hspace{-1mm}&=  D_+ \sum_{j=1}^{m} \left(\begin{matrix}
                \calJ_{i_1*} \\
                \calJ_{i_2*} \\
                \vdots \\
                \calJ_{i_l*} 
            \end{matrix}\right) (-\pz)^j  \frac{d H}{d (\pz ^j x)} (V_{\calI}) \\ 
            &= D_+ \hspace{-1mm}\left(\begin{matrix}
                \calJ_{i_1*} & \ldots & \calJ_{i_1*} \\
                \calJ_{i_2*} & \ldots & \calJ_{i_2*} \\
                \vdots & & \vdots \\
                \calJ_{i_l*} & \ldots & \calJ_{i_l* }
            \end{matrix}\right) \hspace{-2mm}
            \left(\begin{matrix}
                (-\pz)^{1}  \frac{d H}{d (\pz ^{1} x)} (V_{\calI}) \\
                (-\pz)^{2}  \frac{d H}{d (\pz ^{2} x)} (V_{\calI})  \\
                \vdots \\
                (-\pz)^{m_d}  \frac{d H}{d (\pz ^{m_d} x)} (V_{\calI}) 
            \end{matrix}\right).
\end{align*}
As the vector of derivatives of the density function $H$ vanishes for terms not contained in $V_\calI$, we may remove the corresponding columns, which together with \eqref{eq:lifted_Ham} implies  
\begin{equation*}
     \frac{\mathrm{d}}{\mathrm{d}t} \Bar{x} 
     =D_+\hspace{-1mm} \left( \hspace{-0.8mm}\begin{matrix}
                \calJ_{i_1i_1} & \ldots & \calJ_{i_1i_l} \\
                \vdots & \ddots & \vdots \\
                \calJ_{i_li_1} & \ldots & \calJ_{i_li_l}
            \end{matrix}\hspace{-0.8mm}\right)\hspace{-1mm} D_- \delta_{\bar{x}} \bar{\calH} (\bar{x})=\J  \delta_{\bar{x}} \bar{\calH} (\bar{x}).
\end{equation*}
Integration by parts yields $D_{\pm}^{*}=D_{\mp}$, i.e., $\J=-\J^*$.
\end{proof}

\noindent In Proposition \ref{prop:Extension1} we derived from the system \eqref{eq:ham_sys} an augmented Hamiltonian system on the augmented state space with a Hamiltonian not involving differential operators. However, the extended Hamiltonian operator \eqref{eq:ExtendedHamiltonianOperator} is not yet written as matrix differential operator \eqref{eq:Hamiltonian_operator} which is used to define the boundary port variables and the Stokes-Dirac structures in \cite{le2005dirac,villegas2007port,MascvdSc23}. Therefore we derive the decomposition \eqref{eq:Hamiltonian_operator} in the next proposition.

\begin{prop}\label{prop:coeff} 
The differential operator in \eqref{eq:ExtendedHamiltonianOperator} satisfies
    \begin{equation*}
        \J=\sum\limits_{k=0}^{m+2j_l} J_k \pz^k, \qquad \mathrm{with}\quad J_k \coloneqq\sum\limits_{|(u,v,w)|=k} \!\!\!\!Q_{uvw} \in \R^{l\times l}
    \end{equation*}
where for $(u,v,w) \in [0:m] \times \{j_1, \ldots, j_l\}^2$,
\begin{align*}
     (Q_{uvw})_{ij}=\begin{cases}
        (-1)^{v}(P_u)_{i(v),i(w)} \ \ \ &(i,j)=(\operatorname{ind}(v),\operatorname{ind}(w)) \\
        0 \ \ \  & (i,j) \neq (\operatorname{ind}(v),\operatorname{ind}(w))
    \end{cases}
\end{align*}
\end{prop}
\begin{proof}
We rewrite $\J$ as
    \begin{equation*}
        \begin{array}{rcl}
            \J&=& D_+ \left(\begin{matrix}
                \calJ_{i_1i_1} & \ldots & \calJ_{i_1i_l} \\
                \vdots & \ddots & \vdots \\
                \calJ_{i_li_1} & \ldots & \calJ_{i_li_l}
            \end{matrix}\right) D_- \\ [4.5ex]
            &=& D_+ \left( 
            \sum\limits_{j=0}^m \left(\begin{matrix}
                (P_j)_{i_1i_1} & \ldots & (P_j)_{i_1i_l} \\
                \vdots & \ddots & \vdots \\
                (P_j)_{i_li_1} & \ldots & (P_j)_{i_li_l}
            \end{matrix}\right) \pz^j
            \right) \hspace{-1ex}  
            D_- \\[4.5ex]
            &=& \sum\limits_{j=0}^m D_+ \left(\begin{matrix}
                (P_j)_{i_1i_1} & \ldots & (P_j)_{i_1i_l} \\
                \vdots & \ddots & \vdots \\
                (P_j)_{i_li_1} & \ldots & (P_j)_{i_li_l}
            \end{matrix}\right) D_- \pz^j.
        \end{array}
    \end{equation*}
    Since we multiply the coefficient matrix with derivatives of different order from the right and left side, lengthy but straightforward computations yields the result, as~\begin{equation*}
        \begin{array}{rcl}
            \J&=& \sum_{u=0}^m \sum_{v=j_1}^{j_l} \sum_{w=j_1}^{j_l} Q_{uvw} \pz^{u+v+w}=\sum\limits_{k=0}^{m+2j_l} J_k \pz^k.
        \end{array}
    \end{equation*}
\end{proof}
Next we derive the Hamiltonian operators defining the lift of the Hamiltonian formulation for the flexible rod and the Boussinesq equation according to  Proposition~\ref{prop:coeff}.

\begin{ex}
The Hamiltonian operator of the Boussinesq equation is given by 
\begin{align*}
    \left(\begin{matrix}
        0 & \pz\\
        \pz & 0
    \end{matrix}\right)=\left(
    \begin{matrix}
        0 & 0\\
        0 & 0
    \end{matrix}\right)\pz^{0}+\left(
    \begin{matrix}
        0 & 1\\
        1 & 0
    \end{matrix}\right)\pz^{1}\eqqcolon 
    P_{0}\pz^{0}+P_{1}\pz^{1}
\end{align*}
Further, the index set $\calI_{B}$ associated with the Boussinesq equation is $\calI_{B} \coloneqq \{(1,0),(2,0),(1,1)\}$. The augmented state and the lifted Hamiltonian function are defined in \eqref{eq:Boussinesq_state} and \eqref{eq:Boussinesq_Hamiltonian}, respectively. By Proposition \ref{prop:coeff}, we represent the corresponding Hamiltonian operator by $\J_{B}=\sum\limits_{k=0}^{1+2\cdot1}J_k\pz^{k}$ with
{\allowdisplaybreaks
\begin{align*}
    J_0 
            &= \left( \begin{smallmatrix}
                (-1)^0 (P_0)_{11} & 0 & 0 \\ 0 & 0 & 0 \\ 0 & 0 & 0 
            \end{smallmatrix}\right) + \left( \begin{smallmatrix}
                0 & (-1)^0 (P_0)_{12} & 0 \\ 0 & 0 & 0 \\ 0 & 0 & 0 
            \end{smallmatrix}\right) 
            + \left( \begin{smallmatrix}
                0 & 0 & 0 \\ (-1)^0 (P_0)_{21} & 0 & 0 \\ 0 & 0 & 0 
            \end{smallmatrix}\right) + \left( \begin{smallmatrix}
                0 & 0 & 0 \\ 0 & (-1)^0 (P_0)_{22} & 0 \\ 0 & 0 & 0 
            \end{smallmatrix}\right) = 0_3, \\ 
            J_1 
            &= \left( \begin{smallmatrix}
                (-1)^0 (P_1)_{11} & 0 & 0 \\ 0 & 0 & 0 \\ 0 & 0 & 0 
            \end{smallmatrix}\right) + \left( \begin{smallmatrix}
                0 & (-1)^0 (P_1)_{12} & 0 \\ 0 & 0 & 0 \\ 0 & 0 & 0 
            \end{smallmatrix}\right) 
            + \left( \begin{smallmatrix}
                0 & 0 & 0 \\ (-1)^0 (P_1)_{21} & 0 & 0 \\ 0 & 0 & 0 
            \end{smallmatrix}\right) + \left( \begin{smallmatrix}
                0 & 0 & 0 \\ 0 & (-1)^0 (P_1)_{22} & 0 \\ 0 & 0 & 0 
            \end{smallmatrix}\right)
            + 4 \cdot 0_3 \\
            &= \left( \begin{smallmatrix}
                0 & 1 & 0 \\ 0 & 0 & 0 \\ 0 & 0 & 0 
            \end{smallmatrix}\right) + \left( \begin{smallmatrix}
                0 & 0 & 0 \\ 1 & 0 & 0 \\ 0 & 0 & 0 
            \end{smallmatrix}\right)=\left( \begin{smallmatrix}
                0 & 1 & 0 \\ 1 & 0 & 0 \\ 0 & 0 & 0 
            \end{smallmatrix}\right),
\end{align*}
}
and analogously $J_2 = \left( \begin{smallmatrix}
                0 & 0 & 0 \\ 0 & 0 & 1 \\ 0 & -1 & 0 
            \end{smallmatrix}\right)$ and $J_3 = 0_3$. Thus, we obtain the Hamiltonian operator
    \begin{align*}
        0_3\pz^0 +\left( \begin{matrix}
                0 & 1 & 0 \\ 1 & 0 & 0 \\ 0 & 0 & 0 
            \end{matrix}\right) \pz^1 + \left( \begin{matrix}
                0 & 0 & 0 \\ 0 & 0 & 1 \\ 0 & -1 & 0 
            \end{matrix}\right) \pz^2 + 0_3\pz^3 = \left(\begin{matrix}
            0 & \partial_z  & 0 \\ \partial_z & 0 & -\partial_z^2 \\
            0 & \partial_z^2 & 0
        \end{matrix}\right).
    \end{align*}
    \end{ex}

\begin{ex}
Analogously, we can apply Proposition~\ref{prop:coeff} to compute the coefficients of the extended Hamiltonian operator associated to the elastic rod with local elasticity relation as 
\begin{align*}
            \left( \begin{matrix}
                0 & 1 & 0 \\ -1 & 0 & 0 \\ 0 & 0 & 0 
            \end{matrix}\right) \pz^0 + \left( \begin{matrix}
                0 & 0 & 0 \\ 0 & 0 & 1 \\ 0 & 1 & 0 
            \end{matrix}\right)  \pz^1+0_3\pz^2=\left(\begin{matrix}
                0 & 1 & 0 \\ -1 & 0 & \partial_z \\ 0 & \partial_z & 0
            \end{matrix}\right).
        \end{align*}
\end{ex}
    
\section{Dissipative Hamiltonian systems}

\label{sec:diss} 

In this section, we shall extend the previous results to dissipative Hamiltonian systems, introduced in \cite{SCL_16,MacchelliMaschkeGeoplex09} and which are defined with respect to an operator consisting in the sum of an Hamiltonian and a self-adjoint operator. In the sequel, we give a slightly generalized definition by considering a Hamiltonian functional depending on the derivatives of the state variable (its prolongation).
\begin{defn}
Let $\calJ$ a formally skew-adjoint operator of the form \eqref{eq:Hamiltonian_operator} and let $\calG_R^*$ be the formal $L^2$-adjoint of the matrix differential operator 
\begin{align}\label{eq:G_matrix_diff_op}
    \calG_R\coloneqq\sum_{k=0}^{m_g} G_k \partial_z^k
\end{align}
with $G_k \in \R^{n \times d_g}$ for $k \in [0:m_g]$. Further, let $R \in C^\infty([a,b], \R^{d_g \times d_g})$ be such that $\langle e_R, Rf_R\rangle_{L^2([a,b], \R^{d_g})} \geq 0$ for all $(f_R,e_R) \in L^2([a,b], \R^{d_g})$. Then,
\begin{align}\label{eq:diss_sys}
    \dot{x} = (\calJ - \calG_RR\calG_R^*)\delta_x \calH (V_{\calI})
\end{align}
is called a \textit{dissipative Hamiltonian system}.
\end{defn}

It is important to note that the system \eqref{eq:diss_sys} may be expressed by the algebro-differential system
\begin{align}\label{eq:intercon_diss_sys}
    \begin{pmatrix}
        \dot{x}\\
        F
    \end{pmatrix}=
    \begin{pmatrix}
        \calJ & \calG_{R}\\
        -\calG_{R}^{*} & 0
    \end{pmatrix}\begin{pmatrix}
        \delta_{x}\calH(V_{\calI})\\
        \varphi
    \end{pmatrix}
\end{align}
together with the closure relation $\varphi=RF.$ This formulation reveals that there is again a skew-adjoint operator governing the right-hand side. 
For matrix differential operators with constant coefficients it is a Hamiltonian operator, at the heart of the model. This operator induces boundary port variables and a Stokes-Dirac structure \eqref{eq:StokesDiracConstEq}. 

In the sequel, we will embed the dissipative Hamiltonian system \eqref{eq:diss_sys} into the jet space to eliminate the spatial derivatives in the Hamiltonian density and show that this embedding defines again a dissipative Hamiltonian system. 

\begin{prop}\label{prop:Extension_diss}
The extension of \eqref{eq:diss_sys} onto the jet space is again a dissipative 
Hamiltonian system, given by
\begin{equation}\label{eq: explicit_jet_dynamics_diss}
    \frac{\mathrm{d}}{\mathrm{d}t}  
    \bar{x}
    =(\J - \mathbb{G}_R R \mathbb{G}_R^*)  \delta_{\bar{x}} \bar{\calH}(\bar{x}) 
\end{equation}
with $\J$ as in \eqref{eq:ExtendedHamiltonianOperator} and
\begin{equation}\label{eq:Gdef}
     \mathbb{G}_R = D_+ \left(
    \sum_{k=0}^{m_g}
    \left(\begin{matrix}
        (G_k)_{i_1*} \\ \vdots \\ (G_k)_{i_l*}  
    \end{matrix}\right) \partial_z^k \right).
\end{equation}
\end{prop}
\begin{proof}
    The system dynamics \eqref{eq:diss_sys} may be rewritten to
\begin{align*}
    \dot{x} 
    &=\calJ\delta_x \calH (V_{\calI}) - \calG_RR\calG_R^*\delta_x \calH (V_{\calI}).
\end{align*}
Following the proof of Proposition \ref{prop:Extension1} yields that
\begin{align*}
    \frac{\mathrm{d}}{\mathrm{d}t}  
    \bar{x}
    &=\mathbb{J} \delta_{\bar{x}} \bar{\calH}(\bar{x}) - \mathbb{D} \delta_{\bar{x}} \bar{\calH}(\bar{x}),
\end{align*}
with
    $$\mathbb{D}\!=\!D_+ \!\!\left( \begin{matrix}
        (\calG_R R \calG_R^*)_{i_1i_1} & \ldots & (\calG_R R \calG_R^*)_{i_1i_l} \\
        \vdots & \ddots & \vdots \\
        (\calG_R R \calG_R^*)_{i_li_1} & \ldots & (\calG_R R \calG_R^*)_{i_li_l} 
    \end{matrix} \right)\!\!D_- \!=\! \mathbb{G}_R R \mathbb{G}_R^*.$$
\end{proof}

We now provide a representation of $\mathbb{G}_{R}$ in terms of matrix differential operators of the form \eqref{eq:G_matrix_diff_op}. 

\begin{prop}\label{prop:Gcoeff}
    Let $\calG_R\coloneqq\sum_{k=0}^{m_g} G_k \partial_z^k$ with $G_k \in \R^{d_g \times n}$ for $k \in [0:m_g]$ be a matrix differential operator and let $\mathbb{G}_R$ be given by \eqref{eq:Gdef}. Then 
    \begin{align}
        \mathbb{G}_R=\sum_{k=0}^{m_g+j_l} \left(\sum_{|(u,w)|=k} U_{uw} \right) \partial_z^k\eqqcolon\sum_{k=0}^{m_g+j_l} H_k \partial_z^k.
    \end{align}
     where for all $(u,w) \in \{0, \ldots, m_g\} \times \{j_1, \ldots, j_l\}$,
\begin{align*}
    (U_{uw})_{i*}=\begin{cases}
        (G_u)_{i(w)*} \quad & i=\operatorname{ind}(w) \\
        0 \quad & i \neq \operatorname{ind}(w)
    \end{cases}
\end{align*}
\end{prop}
\begin{proof}
The claim follows analogously to Proposition~\ref{prop:coeff} from direct computations.
\end{proof}
Last, we specify the coefficients of the formally skew-symmetric matrix polynomial in \eqref{eq:intercon_diss_sys} defining the Stokes-Dirac structure for the geometric formulation of \eqref{eq:diss_sys}. 
\begin{prop}
consider the operators $\J$ and $\mathbb{G}_{R}$ defined in \eqref{eq:ExtendedHamiltonianOperator} and \eqref{eq:Gdef} and let $\calJ=\sum_{j=0}^{m}J_{j}\partial_{z}^{j}$ and $\calG_{R}=\sum_{k=0}^{m_{g}}G_{k}\partial_{z}^{k}$. Then with $\widehat{m}\coloneqq\operatorname{max}\{m+2j_{l},m_{g}+j_{l})\}$,
\begin{align}\label{eq:extended_ham_op_diss}
    \begin{pmatrix}
        \J & \mathbb{G}_{R}\\
        -\mathbb{G}_{R}^{*} & 0
    \end{pmatrix}=\sum_{k=0}^{\widehat{m}}
    \begin{pmatrix}
        R_{k} & H_{k}\\
        (-1)^{k}H_{k}^{T} & 0
    \end{pmatrix}\partial_{z}^{k}.
\end{align}

\end{prop}
\begin{proof} 
First, note that the formal adjoint of $\mathbb{G}_{R}$ is given by $\mathbb{G}_{R}^{*}=\sum_{k=0}^{m_{g}+j_{l}}H_{k}^{T}(-\partial_{z})^{k}$. Then, the claim follows directly by applying Proposition \ref{prop:coeff} and Proposition \ref{prop:Gcoeff}. 
\end{proof} 

We now conclude this section by an example. 

\begin{ex}[Allen-Cahn equation] 
We recall the port-Hamiltonian formulation of the Allen-Cahn equation as given in \cite{vincent2020port, yaghi2022port, bendimerad2022structure}. Let $\phi$ be the phase variable and consider the \textit{Landau-Ginzburg entropy functional} 
\begin{align*}
    G(\phi)\coloneqq\int_{a}^{b}f(\phi(z,t))+\frac{1}{2}\kappa(z)(\partial_{z}\phi(z,t))^{2}\,\mathrm{d}z.
\end{align*}
Following \cite{bendimerad2022structure}, $f:\R\to\R$ encapsulates the bulk free energy associated with a particular phase state. Further, $\kappa(z)$ is a positive parameter that pertains to the thickness of the interface between two phases. Altogether, $G(\phi)$ symbolizes the free energy of the system. Then, the Allen-Cahn equation is given by 
\begin{align}\label{eq:Allen_Cahn}
    \frac{\mathrm{d}}{\mathrm{d}t}\phi(z,t)=-\Gamma_{\phi}(z)\delta_{\phi}G(\phi),
\end{align}
where $\Gamma_{\phi}:[a,b]\to\R$ denotes the isotropic interface mobility and $\delta_{\phi}G(\phi)$ the variational derivative of $G$ with respect to $\phi$. Note that \eqref{eq:Allen_Cahn} is a dissipative Hamiltonian system with $\calG_{R}=\calG_{R}^{*}=1$,
$R=\Gamma_{\phi}(z)$ and $\delta_{\phi}G(\phi)=\delta_{\phi}f(\phi)-\partial_{z}(\kappa(z)\partial_{z}\phi)$. 
To extend this system onto the jet space, we set $\psi(z,t)=\partial_{z}\phi(z,t)$ such that 
\begin{align*}
\bar{G}(\phi,\psi)\coloneqq\int_{a}^{b}f(\phi(z,t))+\frac{1}{2}\kappa(z)\psi(z,t)^{2}\,\mathrm{d}z.
\end{align*}
The corresponding indices are given by $\{(1,0),(1,1)\}.$ 
Then, by Proposition \ref{prop:Extension_diss} and 
\begin{align*}
    \mathbb{G}_{R}=
    \left(\begin{matrix}
        1 & 0\\
        0 & \partial_{z}
    \end{matrix}\right)
    \left(\begin{matrix}
        1\\1
    \end{matrix}\right)=
    \left(\begin{matrix} 
        1\\
        \partial_{z}
    \end{matrix}\right),\end{align*}
we can reformulate the Allen-Cahn equation as 
\begin{align}\label{eq:extended_AC}
    \frac{\mathrm{d}}{\mathrm{d}t}
    \left( \hspace{-1mm} \begin{matrix} 
        \phi(z,t)\\
        \psi(z,t)
    \end{matrix} \hspace{-1mm} \right) \!=\!-\!
    \left( \hspace{-0.5mm} \begin{matrix} 
        1\\
        \partial_{z}
    \end{matrix} \hspace{-0.5mm} \right) \Gamma_{\phi}(z)
    \left( \hspace{-0.5mm} \begin{matrix} 
        1\\
        \partial_{z}
    \end{matrix} \hspace{-0.5mm} \right)^{*}%
\left( \hspace{-1mm} \begin{matrix} 
    \delta_{\phi}\bar{G}(\phi,\psi)\\
    \delta_{\psi}\bar{G}(\phi,\psi)
\end{matrix} \hspace{-1mm} \right).
\end{align}
Now we are able to express the Allen-Cahn equation using Stokes-Dirac structure combined with a dissipative relation, since the entropy function $\bar{G}$ does not depend on spatial derivatives of the (extended) state $x=(\phi,\psi)$. At first, note that \eqref{eq:extended_AC}
is equivalent to 
\begin{align}\label{eq:extended_AC_equ}
    \left( \hspace{-1mm} \begin{matrix} 
            \tfrac{\mathrm{d}}{\mathrm{d}t}\phi(z,t)\\
            \tfrac{\mathrm{d}}{\mathrm{d}t}\psi(z,t)\\
            F
    \end{matrix} \hspace{-1mm} \right)\hspace{-1mm}  =\hspace{-1mm} 
    \left( \hspace{-0.5mm} \begin{matrix}
        0 & 0 & 1\\
        0 & 0 & \partial_{z}\\
        -1 & \partial_{z} & 0
    \end{matrix} \hspace{-0.5mm}  \right) \hspace{-1mm} 
    \left( \hspace{-1mm} \begin{matrix}
        \delta_{\phi}\bar{G}(\phi,\psi)\\
        \delta_{\psi}\bar{G}(\phi,\psi)\\
        \varphi
    \end{matrix}\hspace{-1mm}  \right) \hspace{-1mm} \eqqcolon\hspace{-0.5mm}  \calJ \left( \hspace{-1mm}  \begin{matrix}
        \delta_{\phi}\bar{G}(\phi,\psi)\\
        \delta_{\psi}\bar{G}(\phi,\psi)\\
        \varphi
    \end{matrix}\hspace{-1mm}  \right)
\end{align}
together with the closure relation $\varphi=-\Gamma_{\phi}(z)F$. To represent the system by geometric structures it is unavoidable to define boundary ports according to the definition given in Section \ref{sec:problem}. Note that the matrix differential operator $\calJ$ is the same as in the case of the elastic rod where we refer to reader to \cite{MascvdSc23} for more details. This is why the boundary ports are given by 
\begin{align}\label{eq:bdd_ports_AC}
    \left(\begin{matrix}
        f_\partial \\ e_\partial 
    \end{matrix}\right)=  \left(\begin{matrix}
        0 & 0 & 1 \\ 0 & 1 & 0
    \end{matrix}\right) \tr  \left(\begin{matrix}
            \delta_\phi \bar{G}(\phi, \psi) \\ \delta_\psi \bar{G}(\phi, \psi)  \\ \varphi
        \end{matrix} \right)
\end{align}
Let $\calB= (L^2([a,b],\R^3) \times \R^{2})^2$. The geometric formulation takes into account the associated Stokes-Dirac structure
\begin{align*}
        \calD_\calJ =
    \left\{ (f,f_{\delta},e, e_{\delta}) \in \calB  \middle| \hspace{-0.5ex}
    \begin{array}{l}
        e \in H^1([a,b], \R^3), f=\calJ e  ,  \\ 
        \begin{bmatrix}
            f_\delta \\ e_\delta
        \end{bmatrix} = \begin{pmatrix}
        0 & 0 & 1 \\ 0 & 1 & 0
    \end{pmatrix} \tr e \end{array} \right\}
\end{align*}
together with the resistive relation
\begin{align*}
        \calR = \left\lbrace (f_R, e_R) \in L^2([a,b], \R )^2 \, \middle| \, e_R = -\Gamma_\phi (z) f_R \right\rbrace.
\end{align*}
Finally, we can describe the dynamics of the Allen-Cahn equation with state by
$(\dot{x}, e_R, f_\delta, \delta_x \bar{G}, f_R, e_\delta) \in \calD_\calJ$ and $(f_R,e_R) \in \calR$. Note that the formally skew-adjoint matrix differential operator $\calJ$ defining the Stokes-Dirac structure of the Allen-Cahn equation is given by \eqref{eq:extended_ham_op_diss}.
\end{ex}

\section{Conclusion}\label{sec:outlook}
In this work we provided an analytical foundation for jet-space extensions of infinite-dimensional conservative and dissipative Hamiltonian systems. More precisely, we showed that jet-space extensions of Hamiltonian systems again yield Hamiltonian dynamics. These jet-space extensions allow for a formulation on a lifted space of extended state variables in which the Hamiltonian density does not include derivatives of the state variable. As a major consequence, this allows for a straightforward derivation of port variables and a corresponding energy-balance.

\section{Acknowledgments}
We are indebted to Friedrich Philipp and Karl Worthmann (TU Ilmenau) for valuable discussions on port-Hamiltonian systems.

\bibliographystyle{plain}
\bibliography{references}

\end{document}